\documentclass[12pt]{article}
\usepackage{amsmath}
\usepackage{amssymb}
\usepackage{graphicx}

\newtheorem{theo}{Theorem}[section]
\newtheorem{lm}{Lemma}[section]
\newtheorem{df}{Definition}[section]

\allowdisplaybreaks

\numberwithin{equation}{section}

\def\R{{\mathbb R}}
\def\Z{{\mathbb Z}}

\def\Pr{{\mathbf P}}

\def\1{{\mathbf 1}}

\def\eps{\varepsilon}
\let\phi=\varphi
\def\qed{\hfill$\Box$}

\title{Percolation for the stable marriage of Poisson and Lebesgue}
\author{M.V.~Freire$^{1}$ \and S.~Popov\thanks{Corresponding author}$^{~,2}$
\and M.~Vachkovskaia$^{1}$}

\begin{document}
\maketitle
{\footnotesize
\noindent
$^1$ Instituto de Matem{\'a}tica, Estat{\'\i}stica e Computa\c{c}\~ao Cient{\'\i}fica,
Universidade de Cam\-pi\-nas,
Caixa Postal 6065, CEP 13083--970, Campinas SP, Brasil.\\
E-mails: mvf@ime.unicamp.br, marinav@ime.unicamp.br

\noindent
$^2$Instituto de Matem{\'a}tica e Estat{\'\i}stica,
Universidade de S{\~a}o Paulo, rua do Mat{\~a}o 1010, CEP 05508--090,
S{\~a}o Paulo SP, Brasil\\
E-mail: popov@ime.usp.br

}

\begin{abstract}
Let~$\Xi$ be the set of points (we call the elements of~$\Xi$ centers) of
Poisson  process in~$\R^d$, $d\geq 2$, with unit intensity. Consider the allocation
of~$\R^d$ to~$\Xi$ which is stable in the sense of Gale-Shapley marriage problem
and in which each center claims a region of volume~$\alpha\leq 1$.
We prove that there is no percolation in the set of claimed sites if~$\alpha$
is small enough, and that, for high dimensions, there is percolation
in the set of claimed sites if~$\alpha<1$ is large enough. \\[.3cm]
\end{abstract}

\section{Introduction and results}
The following model was considered in~\cite{HHP1, HHP2}.
Whenever possible, we keep here the same notation.
The elements of $\R^d$, $d\ge 2$, are called {\it sites}.
We write~$|\cdot |$ for the Euclidean norm and ${\cal L}$ for the
Lebesgue measure in $\R^d$. Let~$\Xi$ be the set of points of
Poisson point process~$\Pi$ in $\R^d$ with intensity $\lambda$
(usually we will assume that $\lambda=1$).
The elements of~$\Xi$ are called {\it centers}.
Let $\alpha\in [0, +\infty]$ be a parameter called the {\it appetite}.
An {\it allocation\/}
of $\R^d$ to~$\Xi$ with appetite~$\alpha$ is a measurable function
$\psi: \R^d\to \Xi\cup \{\infty, \Delta\}$, such that ${\cal L}[\psi^{-1}(\Delta)]=0$,
and ${\cal L}[\psi^{-1}(\xi)]\le \alpha$ for all $\xi\in \Xi$. The set
$\psi^{-1}(\xi)$ is called the {\it territory\/} of the center~$\xi$. We say
that~$\xi$ is {\it sated}\/ if ${\cal L}[\psi^{-1}(\xi)]=\alpha$, and
{\it unsated}\/ otherwise. We say that a site~$x$ is {\it claimed\/} if $\psi(x)\in \Xi$,
and {\it unclaimed}\/ if $\psi(x)=\infty$. 
Here $\psi(x)=\infty$ means that the site~$x$ is unable
to find any center willing to accept it, and $\psi(x)=\Delta$
means that~$x$ is unable to decide between two or more different centers
(intuitively, this means that~$x$
is exactly on the frontier between the territories of different centers
 from~$\Xi$). A ${\cal L}$-null set of sites with
$\psi(x)=\Delta$ is allowed for technical reasons.
% For the purposes of this paper, it is convenient to suppose 
% that $\psi^{-1}(\Delta)$ is empty, i.e., the ties are resolved
% according to some rule.

Stability of an allocation is defined in the following way.
Let~$\xi$ be a center and let~$x$ be a site with $\psi(x)\notin\{\xi, \Delta\}$.
We say that~$x$ {\it desires\/}~$\xi$ if
\[
|x-\xi|<|x-\psi(x)| \text{ or $x$ is unclaimed.}
\]
We say that~$\xi$ {\it covets\/}~$x$, if
\[
|x-\xi|<|x'-\xi| \text{ for some $x'\in\psi^{-1}(\xi)$ or $\xi$ is unsated.}
\]
A site-center pair $(x, \xi)$ is {\it unstable} for the allocation $\psi$
if~$x$ desires~$\xi$ and~$\xi$ covets~$x$.
An allocation is {\it stable\/} if there are no unstable pairs.

\begin{figure}
\centering
\includegraphics[width=8.1cm]{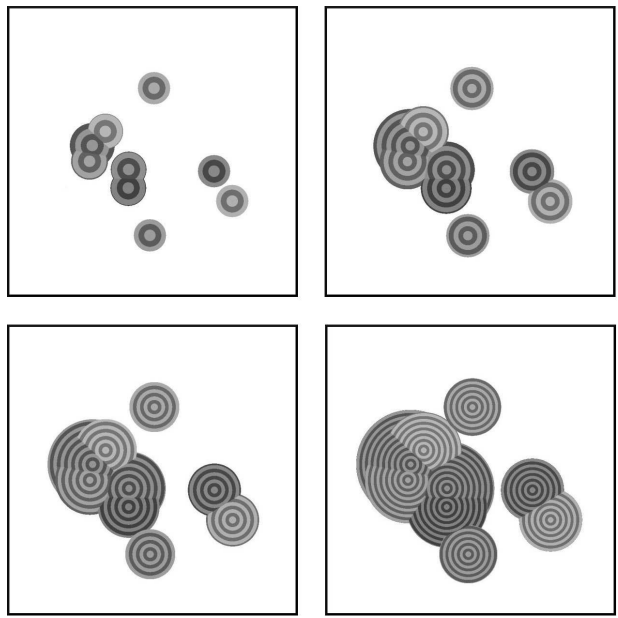}
\caption{Stable allocations for a {\it finite\/} configuration
of centers $\Xi$, and with appetites $\alpha=0.25, 0.45, 0.6, 0.8$.}
\label{alloc}
\end{figure}

The above definition of stable allocation is not constructive. A more 
constructive version can be found in Section 2 of~\cite{HHP1}. 
Informally, the explicit construction of the stable 
allocation can be described as follows. 
For each center, we start growing a ball centered in it. 
All the balls grow simultaneously, at the same linear speed. 
Each center gets the sites captured by its ball, 
unless it is sated or the site was already captured by some other center. 
Remembering  that one picture is worth a thousand words, we refer to 
Figure~\ref{alloc}. Also, it is worth noting that the territory of a particular center 
is not necessarily connected 
(one can imagine the following situation: a center is surrounded 
by several other centers, so the territory 
it gets near itself is not enough, and so it has 
to wait until the neighbouring centers are sated to look for more territory outside).

In~\cite{HHP1}, among other results, it was proved the existence of stable allocation for any set of centers
and any $\alpha\in[0, +\infty]$ and $\Pr$-a.s ${\cal L}$-uniqueness of the stable allocation
in the both following cases:
\begin{itemize}
\item[(i)] $\Xi$ is given by a set of points of an ergodic point process in $\R^d$ or
\item[(ii)] $\Xi$ is finite.
\end{itemize}

Also, it was proved that
\begin{itemize}
\item if $\lambda\alpha<1$ (subcritical regime) then a.s.\ all centers are sated
but there is an infinite volume of unclaimed sites;
\item if $\lambda\alpha=1$ (critical regime) then a.s.\ all centers are sated
and ${\cal L}$-a.a.\ sites are claimed;
\item if $\lambda\alpha>1$ (supercritical regime) then a.s.\ not all centers are sated
but ${\cal L}$-a.a.\ sites are claimed.
\end{itemize}

Denote by
% \[
% {\cal C}'=\{x\in \R^d: \; \psi(x)\in \Xi\} 
% \]
% and let
${\cal C}$ the closure of $\psi^{-1}(\Xi)$. The set ${\cal C}$ 
is the main object of study in this paper; it will be referred to as the 
\emph{set of claimed sites} 
(even though it may contain some $x\in \R^d$ with $\psi(x)=\Delta$).

As shown in~\cite{HHP1},  this model has nice monotonicity properties, both in~$\alpha$ and~$\Xi$ (see Propositions 21 and 22 of~\cite{HHP1}). In this paper, we only need some particular cases of what was proven there, namely,
\begin{itemize}
\item[(i)] if  the sets ${\cal C}_1$ and  ${\cal C}_2$ are 
constructed using the same set of centers
$\Xi$ and different appetites $\alpha_1$ and $\alpha_2$ respectively, 
and $\alpha_1<\alpha_2$, 
then ${\cal C}_1\subset {\cal C}_2$;
\item[(ii)] if the sets ${\cal C}_1$ and  ${\cal C}_2$ are constructed
 using the same appetite $\alpha$ and different sets of centers $\Xi_1$ 
and $\Xi_2$ respectively, 
and $\Xi_1\subset \Xi_2$, then 
${\cal C}_1\subset {\cal C}_2$.
\end{itemize}

In this paper we partially  solve an open problem suggested in~\cite{HHP1}
concerning the percolation of the claimed sites. 
% Note that, without
% restriction of generality, we can suppose that~${\cal C}$
% is a closed set. In fact, we mean that there exists a stable allocation such 
% that~${\cal C}$ is a closed set; this follows from Theorem 5 (see also Theorem 24)
%of~\cite{HHP1}.
\begin{df}
\label{perc}
We say that there is a percolation by claimed sites, 
if there exists an unbounded
% infinite (i.e., which
% is not contained in any finite box) 
connected subset of~${\cal C}$.
% if with positive probability $0$ belongs to an infinite volume connected subset of 
% ${\cal C}$
\end{df}

Due to the monotonicity properties of the model, it is natural to  define
the percolation threshold~$\alpha_p(d)$ in the following way:
\begin{eqnarray*}
 \alpha_p(d) &=& \sup\{\alpha : \\
&&~~~~~\Pr[0 \text{ belongs to an unbounded connected 
subset of } {\cal C}\\
 &&~~~~~~~~~~~~~\text{in the $d$-dimensional model with appetite }\alpha]=0\}.
% \text{ there is no percolation by claimed sites}\\
%  &&~~~~~~~~~~~~~\text{in the $d$-dimensional model with appetite }\alpha\}.
\end{eqnarray*}

\begin{figure}
\centering
\includegraphics[width=12.3cm]{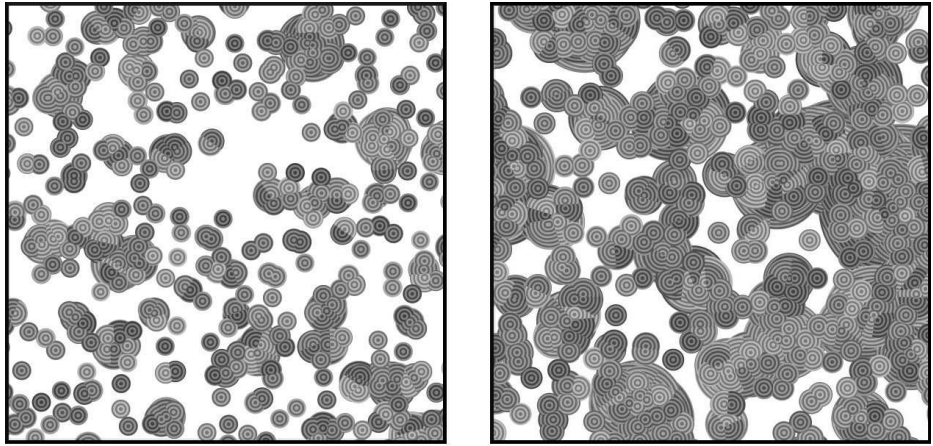}
\caption{On the left image one can see a realization of the model
with $\alpha=0.6$ (which seems to correspond to the non-percolating phase), 
on the right image, $\alpha=0.8$ was used 
(which seems to correspond to the percolating phase).}
\label{06e08}
\end{figure}

On Figure~\ref{06e08} one can see two configurations (inside a box $20\times 20$, 
with $\lambda=1$) in (presumably) non-percolating and percolating phases.

\begin{figure}
\centering
\includegraphics[width=12.3cm]{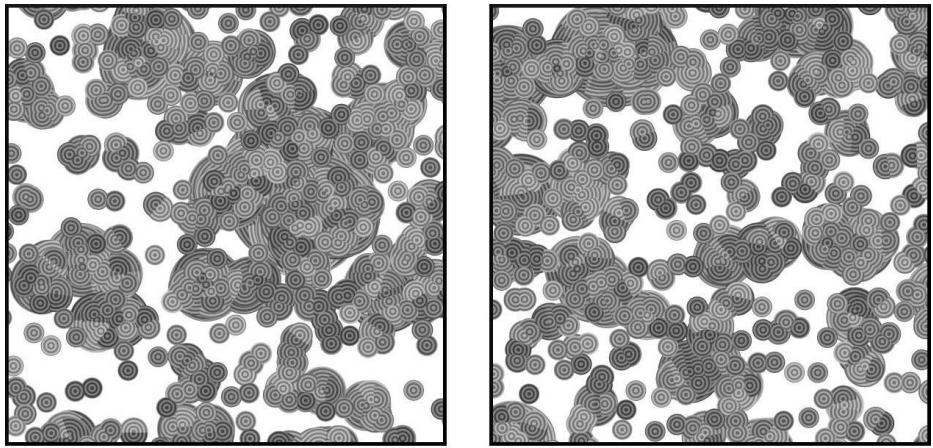}
\caption{Near the percolation threshold (two realizations with $\alpha=0.7$): on the left image, crossings from left to right and from top to bottom do not exist; 
on the right image, 
there are crossings from left to right and from top to bottom.}
\label{07}
\end{figure}

\begin{theo}
\label{t_main}
\begin{itemize}
\item[(i)] For any dimension~$d\geq 2$ we have that $\alpha_p(d)>0$,
that is, if the appetite~$\alpha$ is small enough, then a.s.\
there is no percolation by claimed sites. 
\item[(ii)] Also, if~$\alpha$ is small enough and ~$d\geq 2$, 
then there exists percolation by unclaimed sites 
(i.e., a.s.\ there is an unbounded connected component in $\R^d\setminus {\cal C}$). 
\end{itemize} 
\end{theo}

% 
% \begin{rmk}
% \label{unclaimed}
%  From the proof of Theorem~\ref{t_main} it is not difficult to observe that if
% $\alpha$ is small enough, then there exists percolation by unclaimed sites. 
% \end{rmk}

Since in the model with $\alpha=1$ almost all the sites are claimed, that is, 
${\cal C}=\R^d$, it
is clear that $\alpha_p(d)\leq 1$. 
The next result implies that if the dimension is sufficiently high,
then $\alpha_p(d)<1$ (and even that $\alpha_p(d) \lesssim 2^{-d}$, as $d\to\infty$).
\begin{theo}
\label{t_limsup}
We have
\begin{equation}
\label{limsup_alpha}
\limsup_{d\to\infty}\alpha_p(d)2^{d} \leq 1.
\end{equation}
\end{theo}

%
%Also, we prove the following
%\begin{theo}
%\label{main1}
%Place here results for $d=2$...
%\end{theo}

%Place here simulation results for $\alpha_p$ in the case $d=2$.

%Define the set of unclaimed sites ${\cal U}$ as the closure of $\R^d\setminus {\cal C}$.
%Make a remark about percolation by unclaimed sites (if $\alpha$ is small
%then there is percolation by unclaimed sites).

Simulations suggest that $\alpha_p(2)$ is around $0.7$ (see Figure~\ref{07}).
Note, however, that proving that  $\alpha_p(2)<1$ (as well as
$\alpha_p(d)<1$ for small $d$) is still an open problem.

\section{Proofs}
Since the proof of Theorem~\ref{t_limsup} is much simpler, let us begin by

\medskip
\noindent
{\it Proof of Theorem~\ref{t_limsup}}.
Note that if we rescale the space by factor~$b$ (that is, apply a homothetic map
$x \mapsto bx$), then we obtain the model with the intensity of the Poisson process 
being $\lambda/b^d$ and the appetite~$\alpha b^d$. The geometric properties of the 
allocation  do not change under this transformation, 
and the product of intensity and appetite does not change either.
In particular, this shows that the percolation properties of the model only depend on the product~$\lambda\alpha$.

Let~$\pi_d$ be 
volume of the unit ball in~$\R^d$. Since the volume of the ball of
radius $(\alpha/\pi_d)^{1/d}$ is~$\alpha$, any site 
which is at most $(\alpha/\pi_d)^{1/d}$ far away
 from some center will belong to~${\cal C}$. Indeed, the centers 
want territory of volume~$\alpha$ as close as possible, 
so, for any center, any site~$x$ in the ball of volume~$\alpha$ centered there will be claimed, either by this center, or by another one
(or it may happen that $\psi(x)=\Delta$ so that~$x$ is disputed by two or
more centers, but in this case $x\in {\cal C}$ anyway). 
So, the set of claimed sites~${\cal C}$
dominates the Poisson Boolean model with rate $\lambda=1$ and
radius $(\alpha/\pi_d)^{1/d}$. By the above rescaling argument, that model
is equivalent to the Poisson Boolean model with $\lambda=\alpha \pi_d^{-1} 2^d$
and radius~$1/2$. Let~$\lambda_{cr}(d)$ be the critical rate for the
percolation in the Poisson Boolean model with radius~$1/2$.
Now it is straightforward to obtain that Theorem~\ref{t_limsup}
is a consequence of the following result of~\cite{P}:
\[
 \lim_{d\to\infty} \pi_d \lambda_{cr}(d) = 1.
\]
\qed

\medskip

To prove Theorem~\ref{t_main}, we need some preparations. 

Let us from now on fix~$\alpha=1$ and vary~$\lambda$, instead of fixing $\lambda=1$
and varying $\alpha$ (the rescaling argument in the beginning of the proof of 
Theorem~\ref{t_limsup} allows us to do this). 

The idea of the proof of the part~(i) of Theorem~\ref{t_main}
can be described as follows:
\begin{enumerate}
\item We define a ``discrete'' (i.e., made of cubes of size 1) 
dependent percolation model, and prove (Lemma~\ref{l_terr_cub})
that it dominates the original model, so that it is enough to prove the absence
of percolation in this discrete model.
\item The important properties of the discrete model are provided 
by~(\ref{oc_R1}) and Lemma~\ref{l_indep}.
\item Then, in Definition~\ref{df_passable}, we define the notion of
\emph{passable} level-$m$ cube (a level-$m$ cube is a cube of size~$m$,
see~(\ref{level-m_cube}) below), and we show (Lemma~\ref{l_ind_passable})
that the cubes are passable or not independently if they are far enough
from each other.
\item Using that independence, we prove (Lemmas~\ref{l_An} and~\ref{l_A_n_pass})
that the probability that a bigger cube is passable can be bounded from
above in terms of the probability that a smaller cube is passable.
\item This allows us to prove that, for small enough~$\lambda$,
the probability that a cube is passable tends to 0 as the size 
of the cube goes to infinity. With a little more work, this implies the
absence of percolation.
\end{enumerate}

For $m\ge 1$ and $i=(i^{(1)}, \ldots, i^{(d)})\in \Z^d$, define 
the {\it level}-$m$ cube  $K_i^{m}$ associated with $i$ by
\begin{equation}
\label{level-m_cube}
K_i^{m}=\Big\{x=(x^{(1)}, \ldots, x^{(d)})\in\R^d: -\frac{m}{2}\le  x^{(l)}-mi^{(l)}
              \le \frac{m}{2}\Big\}.
\end{equation}
Note that, in the above definition, 
the quantity~$m$ is not necessarily integer (although it is 
convenient to think about it as such). 
Note also that the union of all level-$m$ cubes 
is $\R^d$, and the intersection of any two 
distinct level-$m$ cubes is either empty, or has zero Lebesgue measure.
We say that two cubes are connected if they have at least one point in common.
Denote by~$\zeta^{(i)}$ the number of centers
in~$K_i^1$, i.e., the cardinality of the set $K_i^1 \cap \Xi$.
At this point we need to introduce more notations. First, we define
the distance between two sets $A,B\subset \R^d$ in a usual way:
\[
\rho(A,B)=\inf_{x\in A, y\in B}|x-y|.
\]
Then, for any $r\ge 0$,  we define a discrete ball $B_i(r)$ by
\[
B_i(r)=\bigcup_{j\in \cal J}K_j^1,\quad\text{ if  }r>0
\]
where
\[
{\cal J}=\{j\in \Z^d: \rho(K_i^1,K_j^1)\le r\},
\]
and
$B_i(0):=\emptyset$. We use the notation  $\lceil x \rceil$ 
 for the smallest integer greater than or equal to~$x$ 
and~$\lfloor x \rfloor$ for the integer part of~$x$.
For each $i\in \Z^d$, define also the random variable (note that $\Z^d\cap K_j^1=j$)
\begin{equation}
\label{def_Ri}
R_i=\inf \Big\{r>0: \sum_{j\in \Z^d\cap B_i(\beta_d r)}\zeta^{(j)}\le \pi_d r^d\Big\},
           \text{ if }\zeta^{(i)}>0
\end{equation}
(here we use the convention $\inf\emptyset=+\infty$),
where, as before,  $\pi_d$ is the volume of the ball with radius~$1$ in~$\R^d$,
\begin{equation}
\label{def_beta_d}
\beta_d=\lceil 3+2\sqrt{d}\pi_d^{1/d}\rceil,
\end{equation}
and $R_i:=0$, if $\zeta^{(i)}$=0.

 We have the following
\begin{lm}
\label{l_terr_cub}
The territories from all centers in $K_i^1$ are contained in
$B_i(R_i)$.
\end{lm}

\medskip
\noindent
{\it Proof of Lemma~\ref{l_terr_cub}}.\/
% As  $ \sum_{j\in B_i(\beta_d R_i)}\zeta^{(j)}\le\pi_d R_i^d$,
% we have that  ${\cal L}(B_i(R_i))$ (i.e., the volume of $B_i(R_i)$, note that 
% ${\cal L}(B_i(R_i))\ge \pi_d\R^d$) is sufficient for 
% the centers in $K_i^1$,
% and even for the centers from $B_i(\beta_d R_i)$.
%
Let us first show that
% \begin{equation}
% \label{stab2.1}
% \max_{x\in K_i^1, y\in B_i(R_i)}|x-y| < \rho(B_i(R_i), \R^d\setminus B_i(\beta_d R_i)). 
% \end{equation}
\begin{equation}
\label{stab2.1}
 \rho(B_i(R_i), \R^d\setminus B_i(\beta_d R_i)) > R_i. 
\end{equation}
Indeed, since $\rho (K_i^1,B_i(R_i))\le R_i$ and the level-$1$ cubes 
have side $1$ and thus diameter~$\sqrt{d}$, we have
\begin{equation}
\label{dist_betad}
\max_{x\in K_i^1, y\in B_i(R_i)} |x-y| \le R_i+2\sqrt{d}.
\end{equation}
For $\rho(B_i(R_i), \R^d\setminus B_i(\beta_d R_i))$, using~(\ref{dist_betad})
and the fact that $\rho (K_i^1,B_i(\beta_d R_i))\le \beta_d R_i$, we obtain 
\begin{equation}
\label{noch_einmal}
\rho(B_i(R_i), \R^d\setminus B_i(\beta_d R_i))\ge \beta_d R_i-(R_i+2\sqrt{d}).
\end{equation}
Finally, note that if $R_i>0$, then there is at least one center in 
$B_i(\beta_d R_i)$, and thus $\pi_d R_i^d\ge 1$, so $R_i\ge \pi_d^{-1/d}$.
So, from~(\ref{dist_betad}) and~(\ref{noch_einmal}), we get that 
if $\beta_d>2+2\sqrt{d}\pi_d^{1/d}$ (by~(\ref{def_beta_d}), this is indeed
the case), then~(\ref{stab2.1}) holds.
% \[
% \max_{x\in K_i^1, y\in B_i(R_i)}|x-y| < \rho(B_i(R_i), \R^d\setminus B_i(\beta_d R_i)). 
% \]

% By~(\ref{stab2.1}), any site in  $B_i(R_i)$ is closer to any center in $K_i^1$ 
% than any site in $\R^d\setminus B_i(\beta_d R_i)$.
% Since the allocation  is stable, this means that
% the centers from $\R^d\setminus B_i(\beta_d R_i)$ will be allowed to have
% territory in $B_i(R_i)$ only when  all centers from $K_i^1$ are sated.
%To prove this in a formal way, 
Now,
suppose that there exist $\xi\in K^1_i$ and $x\in \R^d$ such that
$\psi(x)=\xi$ and $|x-\xi| > R_i$. One can choose a small enough~$\eps$
such that $|x-\xi| > R_i + \eps$ and any site~$z$ with $\rho(z,K^1_i)\leq R_i + \eps$
belongs to~$B_i(R_i)$ (this is possible since $B_i(R_i)$ is a compact set, and
for any~$z'$ from the boundary of $B_i(R_i)$ it holds that $\rho(z',K^1_i) > R_i$,
otherwise the next level-1 cube would be included in~$B_i(R_i)$ too).
There exists~$y$ such that $|y-\xi| \leq R_i+\eps$ (and so $y\in B_i(R_i)$) and
\begin{itemize}
 \item either $\psi(y)=\xi'$ for some $\xi'\in \R^d\setminus B_i(\beta_d R_i)$,
 \item or~$y$ is unclaimed.
\end{itemize}
This is because, by~(\ref{def_Ri}), the number of centers in $B_i(\beta_d R_i)$
is at most $\pi_dR_i^d$, and each one of them wants to claim a territory of volume~1,
but ${\cal L}(\{z\in \R^d : |z-\xi| \leq R_i+\eps\}) > \pi_dR_i^d$.
Now, let us show that $(y,\xi)$ is an unstable pair. Indeed,
\begin{itemize}
 \item $y$ desires~$\xi$, because, by~(\ref{stab2.1}), we have $|y-\xi|<|y-\psi(y)|$, and
 \item $\xi$ covets~$y$, because $|y-\xi|<|x-\xi|$.
\end{itemize}
% In other words, if we suppose that there a site in $\R^d\setminus B_i(R_i)$ which 
% is allocated to 
% a center in $K_i^1$, it means that a center from $\R^d\setminus B_i(\beta_d R_i)$
% has territory in  $B_i(R_i)$ (since the volume of $B_i(R_i)$
% is sufficient for the centers in $B_i(\beta_d R_i)$), and in this case the allocation is unstable
% by~(\ref{stab2.1}).
Thus, the centers from $K_i^1$ will be sated with territory inside
$B_i(R_i)$ and Lemma~\ref{l_terr_cub} is proved.
\qed

\medskip

Lemma~\ref{l_terr_cub} allows us to majorize the original model
by the following (dependent) percolation model:
given the set $\Xi$ of points of Poisson process,
for every $K_i^1$ we paint all the level-$1$ cubes in $B_i(R_i)$ and denote by ${\mathfrak C}$ the set of painted sites. That is, we define
\[
 {\mathfrak C}=\bigcup_{i\in \Z^d}B_i(R_i).
\]

At this point it is important to observe that, by Lemma~\ref{l_terr_cub},  
it holds that ${\cal C}\subset {\mathfrak C}$. Thus, to prove the first part of 
Theorem~\ref{t_main},
it is sufficient to prove the absence of the infinite cluster in  ${\mathfrak C}$
for small $\lambda$.

Let us recall  Chernoff's bound for  Poisson random variable $Z$ with 
parameter $\lambda$:
\begin{equation}
\Pr[Z>a]\le e^{-\lambda g(\lambda/a)},
\end{equation}
where $g(x)=[x-1-\log x]/x$ (note that $g(x)\to +\infty$, as $x\to 0$).
Since ${\cal L}(B_i(\beta_d a))\le (2\beta_d a+3)^d$, by~(\ref{def_Ri}), we have
\begin{eqnarray}
\Pr[R_i>a]&\le & \Pr\Big[ \sum_{j\in B_i(\beta_d a)}\zeta^{(j)}> \pi_d a^d\Big]\nonumber \\
&\le& \exp\Big\{ -\lambda (2\beta_d a+3)^d g\Big(\frac{\lambda (2 \beta_d a+3)^d}{\pi_d a^d}\Big)\Big\}\nonumber\\
&\le& e^{-c(\lambda) a^d}\label{oc_R1}
\end{eqnarray}
where $c(\lambda)\to +\infty$, as $\lambda\to 0$ (a similar argument can be found in the
proof of Proposition~11 from~\cite{HHP2}).

The following simple fact is important for the proof of Theorem~\ref{t_main}:
\begin{lm}
\label{l_indep}
To determine whether the event $ \{R_i\le a\}$ occurs, we only have
to look at the configuration of the centers inside $B_i(\beta_d a)$).
%  \item[(ii)] If $|i-j|>2\beta_d a+3\sqrt{d}$, then the events
% $\{R_i\le a\}$ and $\{R_j\le a\}$ are independent.
% \end{itemize}
\end{lm}

\noindent
{\it Proof of Lemma~\ref{l_indep}.}
This follows immediately from the definition
of~$R_i$ (see~(\ref{def_Ri})). 
\qed

\medskip

% It is important to observe that, due to the definition
% of~$R_i$ (see~(\ref{def_Ri})), the events
% \begin{equation}
% \label{indep}
% \text{$\{R_i\le a\}$ and $\{R_j\le a\}$ are independent, if $|i-j|>2\beta_d a+3\sqrt{d}$}
% \end{equation}
% (indeed, to determine whether the event $ \{R_i\le a\}$ occurs, we only have
% to look at the configuration of the centers inside $B_i(\beta_d a)$).

% We prove now that it is possible to choose $\lambda$
% small enough, so that for all~$n$
% the probability of a  level-$n$ cube to have a crossing by painted level-$1$ cubes will
% be  arbitrarily close to $0$. 
% Clearly, this is enough to guarantee the absence of percolation.

Consider a bounded set $W\subset \R^d$ and 
let $\Xi_W=\Xi\cap W$ (since~$W$ is bounded, 
$\Xi_W$ is a finite set a.s.). As noted above, there exists an a.s.\ unique
stable allocation corresponding to the set of centers~$\Xi_W$. 
We can then construct the set of 
painted sites ${\mathfrak C}|_W$ corresponding to this stable allocation  
analogously to the construction of ${\mathfrak C}$. Namely,
first, we define the random variables $\zeta^{(i)}_W$ as the cardinality of the
set $\Xi\cap K_i^1\cap W$. Then, we define~$R_i^W$ analogously to~(\ref{def_Ri})
(only changing $\zeta^{(\cdot)}$ to~$\zeta^{(\cdot)}_W$), and then we
let ${\mathfrak C}|_W = \cup_{i\in \Z^d} B_i(R_i^W)$.
 From~(\ref{def_Ri}) it is straightforward to obtain that
% Definition of $R_i$ and monotonicity properties of 
% the model imply monotonicity of the set of painted 
% sites in $\Xi$. So, 
${\mathfrak C}|_W\subset{\mathfrak C}$ for any $W\subset \R^d$.
% ;
% on the other hand, from~Lemma~\ref{l_terr_cub} we still obtain that
% ${\cal C}|_W \subset {\mathfrak C}|_W$, where ${\cal C}|_W$ is the set
% of claimed sites for the stable allocation constructed using the centers
% in $\Xi\cap W$.

Let~$K^m_j$ be a level-$m$ cube and define
\[
{\cal A}(K^m_j)=\bigcup_{i: K^m_i\cap K^m_j\ne \emptyset} K^m_i
\]
(so, ${\cal A}(K^m_j)$ is the union of $K^m_j$ with the $3^d-1$
neighbouring level-$m$ cubes.)
We use here some ideas typical for multiscale (fractal) percolation models,
see e.g.~\cite{MPV}.
First, we define the notion of {\it passable\/} cubes.
\begin{df}
\label{df_passable}
A level-$m$ cube $K^m_j$ is \emph{passable}, if
\begin{itemize}
\item[(i)] the set $K^m_j$ intersects a connected
component with diameter at least~$m/2$ of ${\mathfrak C}|_{{\cal A}(K^m_j)}$, and
\item[(ii)] for any $i\in {\cal A}(K^m_j)\cap\Z^d$
we have $R_i<\frac{m}{6(\beta_d+1)}$.
\end{itemize}
\end{df}
Denote by $\|\cdot\|_\infty$ the maximum norm in $\Z^d$ and in~$\R^d$. 
The key observation is that the event ``the level-$m$ cube is passable''
only depends on what happens in finitely many level-$m$ cubes around it.
More precisely:
\begin{lm}
\label{l_ind_passable}
Suppose that $m > 6$ and $\|i-j\|_\infty\geq 5$.
Then the events
$\{K^m_i \text{ is passable}\}$ and $\{K^m_j \text{ is passable}\}$
are independent.
\end{lm}

\noindent
{\it Proof of Lemma~\ref{l_ind_passable}.}
Consider $\ell_1\in {\cal A}(K^m_i)\cap \Z^d$ and $\ell_2\in {\cal A}(K^m_j)\cap \Z^d$.
By Lemma~\ref{l_indep}, the event $\{R_{\ell_k}<\frac{m}{6(\beta_d+1)}\}$
only depends on what happens inside 
$B_{\ell_k}(\frac{\beta_d m}{6(\beta_d+1)})$,
$k=1,2$. Note that $B_{\ell_k}(\frac{\beta_d m}{6(\beta_d+1)})\subset B_{l_k}(m/6)$.
It is then straightforward to check that, if $m > 6$
and $\|i-j\|_\infty\geq 5$, for all such
$\ell_1,\ell_2$ it holds that $B_{\ell_1}(m/6) \cap 
B_{\ell_2}(m/6) = \emptyset$, which concludes the proof
of Lemma~\ref{l_ind_passable}.
% %This follows from Lemma~\ref{l_indep}~(i)
\qed

\medskip

% The key observation is
% that, by Lemma~\ref{l_indep}~(ii), if $\|i-j\|_\infty\geq 2 \beta_d+1$, then the events
% $\{K^m_i \text{ is passable}\}$ and $\{K^m_j \text{ is passable}\}$
% are independent ({\bf make it lemma?}) (if $m$ is so large that $m/3>3\sqrt{d}$).

Denote $p_m:=\Pr[\text{$K^m_0$ is passable}]$. Next, our goal is
to show that if~$\lambda$ is small enough, then $p_m\to 0$ as~$m\to\infty$.

Consider the event
\begin{eqnarray*}
A_n&=&\Big\{\text{in ${\cal A}(K^n_0)$ there exists a connected component of 
           diameter}\\
&&~~~~~~~~~~~~~~~~~~~~~~~~\text{at least
$\displaystyle\frac{n}{2}$ of passable level-$(3\log n)$ cubes}\Big\}.
\end{eqnarray*}

\begin{lm}
\label{l_An}
We have, for $n>6$,
\begin{equation}
\label{sy}
\Pr[A_n]\le %\Big(\frac{n}{3\log n}\Big)^d (4\beta_d+3)^{dk_0} p_{3\log n}^{k_0} =
\Big(\frac{n}{\log n}\Big)^d (11^d p_{3\log n})^{k_0},
\end{equation}
where
\[
 k_0=\Big\lfloor\frac{n}{30\sqrt{d}\log n}\Big\rfloor-1.
\]
\end{lm}

\medskip
\noindent
{\it Proof of Lemma~\ref{l_An}\/}.
Since the diameter of a level-$m$ cube is $m\sqrt{d}$, on the event~$A_n$,
there exist~$m'\in \Z_+$, $i_1,\ldots,i_{m'}\in\Z^d$ such that
\begin{itemize}
 \item $K_{i_j}^{3 \log n} \subset {\cal A}(K^n_0)$ for all $j=1,\ldots,m'$,
 \item $\|i_j-i_{j-1}\|_{\infty}=1$, for all $j=2,\ldots,m'$,
 \item and $\|i_0-i_{m'}\|_{\infty}\geq \frac{n}{6\sqrt{d}\log n}$.
\end{itemize}
Then, define $\tau(1):=1$, and
\[
 \tau(j) = \max\{ \ell > \tau(j-1) : \|i_\ell-i_{\tau(j-1)}\|_{\infty} = 5 \}
\]
for $j=2,\ldots,k_0$ (indeed, since $5k_0 < \frac{n}{6\sqrt{d}\log n}$,
we have that $\tau(k_0)\leq m'$). 
Then, the
% inside ${\cal A}(K^n_0)$ there is a connected component
% of passable level-$(3\log n)$ cubes of 
% length $\lfloor\frac{n}{6\sqrt{d}\log n}\rfloor$.
% If the event~$A_n$ occurs, then there exists
% a 
collection of level-$(3\log n)$ cubes
 $\gamma=(K^{3\log n}_{i_{\tau(1)}}, \ldots, K^{3\log n}_{i_{\tau(k_0)}})$
% such that
has the following properties:
$\gamma\subset {\cal A}(K^n_0)$, for $j=1,\ldots,k_0$ the cubes
 $K^{3\log n}_{i_{\tau(j)}}$ are passable, $\|i_{\tau(j)}-i_{\tau(j-1)}\|_\infty=5$
and $\|i_{\tau(j')}-i_{\tau(j)}\|_\infty\ge 5$, for all $j\ne j'$.
Intuitively, this collection corresponds to a ``path'' by
passable cubes inside ${\cal A}(K^n_0)$; however, neighbouring
elements of this path are not really neighbours, but they are separated enough
to make them independent. 
The number of collections with such properties is 
at most $\big(\frac{n}{\log n}\big)^d 11^{dk_0}$
(there are at most $\big(\frac{n}{\log n}\big)^d$ possibilities 
to choose the first cube in the collection,
and then at each step there are at most $11^d-9^d<11^d$ possibilities
to choose the next one).
For a fixed~$\gamma$, by Lemma~\ref{l_ind_passable}, 
the probability that all the cubes $K^{3\log n}_{i_j}$ in the collection~$\gamma$
are passable is at most $p_{3\log n}^{k_0}$, 
and so~(\ref{sy}) holds.
% by independence (see the observation just
% after Definition~\ref{df_passable}).
\qed
% Thus,
% \[
% \Pr[A_n]\le %\Big(\frac{n}{3\log n}\Big)^d (4\beta_d+3)^{dk_0} p_{3\log n}^{k_0} =
% \Big(\frac{n}{\log n}\Big)^d ((4\beta_d+3)^d p_{3\log n})^{k_0},
% \]
% which is small, if $p_{3\log n}<(4\beta_d+3)^{-d}$ and $n$ (and thus $k_0$) is large.

\begin{lm}
\label{l_A_n_pass}
Suppose that, for some $n>e^{4\sqrt{d}}$ the cube~$K^n_0$ is passable
% .
% Suppose also that $K^n_0$ intersects with a connected component with a diameter
% at least~$n/2$ of painted level-$1$ cubes, 
and the following event occurs:
\begin{equation}
\label{Rlog}
\Big\{\text{for all } i\in {\cal A}(K^n_0) \cap \Z^d
  \text{ it holds that }R_i<\frac{\log n}{2(\beta_d+1)}\Big\}.
\end{equation}
Then, any level-$(3\log n)$ cube in 
${\cal A}(K^n_0)$ intersecting with a connected component
 of ${\mathfrak C}|_{{\cal A}(K^n_0)}$ with a diameter
at least~$n/2$ of painted level-$1$ cubes (cf. Definition~\ref{df_passable}~(i)), 
and such that the distance from it to~$K^n_0$ is 
at most~$n/2$, is passable, and, in particular, 
the event~$A_n$ occurs.
\end{lm}

\medskip
\noindent
{\it Proof of Lemma~\ref{l_A_n_pass}}.\/
Consider any level-$(3\log n)$ cube with the above properties, say $K^{3\log n}_j$. 
As $\rho( K^n_0, K^{3\log n}_j)\le n/2$, we have 
$ {\cal A}(K^{3\log n}_j)\subset {\cal A}(K^n_0)$ and thus 
for all $x \in {\cal A}(K^{3\log n}_j)$ it holds that 
$R_x<\frac{\log n}{2(\beta_d+1)}$. 
That is, the second condition in Definition~\ref{df_passable} is satisfied.
Let 
\[
{\mathfrak K}( K^{3\log n}_j)= \big\{x\in \R^d : \inf_{y\in K^{3\log n}_j}
             \|x-y\|_\infty \le 2\log n \big\}.
%\bigcup_{l:\|j-l\|_\infty\le 2\log n} K^{3\log n}_l.
%\{K^{3\log n}_l: \|j-l\|_\infty\le 2\log n \}.
\]
Abbreviate $r_1:=\frac{\log n}{2(\beta_d+1)}$, and
consider some level-1 cube $K^1_l\subset {\mathfrak K}( K^{3\log n}_j)$
such that $K^1_l\subset {\mathfrak C}|_{{\cal A}(K^n_0)}$.
On the event~(\ref{Rlog}) 
this means that there exists $i\in {\cal A}(K^n_0)\cap \Z^d$
such that $K^1_l\subset B_i^{{\cal A}(K^n_0)}(r_1)$.
By Lemma~\ref{l_indep}, 
the event $\{R_i \leq r_1\}$ only depends on the configuration 
inside $B_i(\beta_d r_1)$. Since 
$r_1+\beta_d r_1+2\sqrt{d} = \frac{\log n}{2}+2\sqrt{d} < \log n$
(we supposed that $\log n > 4\sqrt{d}$),
we obtain that $B_i(\beta_d r_1)$ is fully inside ${\cal A}(K^{3\log n}_j)$.
So, 
%Then, by~(\ref{Rlog}), there exists~$j$ such that
$K^1_l\subset B_i \big(R_i^{{\cal A}(K^{3\log n}_j)}\big)$,
and, consequently,
 $ {\mathfrak K}( K^{3\log n}_j)$ intersects 
with a connected component of diameter
at least $2\log n$ of level-$1$ cubes from ${\mathfrak C}|_{{\cal A}(K^{3\log n}_j)}$.
This implies that  $K^{3\log n}_j$ is passable. 
% 
% By~(\ref{Rlog}), all painted sites in $ {\mathfrak K}( K^{3\log n}_j)$ 
% we painted by some $x \in {\cal A}(K^{3\log n}_j)$ (as we suppose 
% that $K^n_0$ is passable,  we already consider  
% the configuration restricted  on ${\cal A}(K^n_0)$).
% So, $ {\mathfrak K}( K^{3\log n}_j)$ intersects with a connected component of diameter
% at least $2\log n$ of level-$1$ cubes painted by some $x \in {\cal A}(K^{3\log n}_j)$, 
% this implies that  $K^{3\log n}_j$ is passable. 
\qed

\medskip

Now we are ready to prove~Theorem~\ref{t_main}.

\medskip
\noindent
{\it Proof of Theorem~\ref{t_main}}. 
Using first Lemma~\ref{l_A_n_pass}, and then~(\ref{oc_R1})
together with Lemma~\ref{l_An}, we obtain that
\begin{eqnarray}
p_n&=&\Pr[\text{$K^n_0$ is passable}]\nonumber\\
&\le&\Pr[A_n]
+\Pr\Big[\text{there exists } i\in 
       {\cal A}(K^n)\cap\Z^d: R_i\ge\frac{\log n}{2(\beta_d+1)}\Big]\nonumber\\
&\le&
\label{pass}
\Big(\frac{n}{\log n}\Big)^d (11^d p_{3\log n})^{k_0}
  + (3n)^d e^{-c'(\lambda) \log^d n},
\end{eqnarray}
where $c'(\lambda) = 2^{-d}(\beta_d+1)^{-d}c(\lambda)$.
% Denote
% \[
% f(n, \eps):= \Big(\frac{n}{\log n}\Big)^d ((4\beta_d+3)^d\eps)^{k_0}
%   + (3n)^d e^{-c(\lambda) \log^d n}
% \]
Abbreviate $\eps_d=11^{-d}/2$. Choose a large enough~$m_0$ such that
\begin{eqnarray*}
 m &> & e^{4\sqrt{d}},\\
m^{1/2}&\leq &\Big\lfloor\frac{m}{30\sqrt{d}\log m}\Big\rfloor-1,\\
\frac{m}{\log m} &\le & 3m,\\
(11^d\eps_d)^{m^{1/2}} &\le & e^{-\log^d m},\\
\eps_d & > &(3m)^d \Big( (11^d\eps_d)^{m^{1/2}}+e^{-\log^d m}\Big)
\end{eqnarray*}
for all $m\ge m_0$ (note that in fact $11^d\eps_d = 1/2$ and that
$e^{4\sqrt{d}}>6$).
Choose a small enough~$\lambda$ in such a way that $c'(\lambda)\ge 1$ 
and also that $e^{-\lambda(3m_0)^d}>1-\eps_d$.
Note that the last condition on~$\lambda$ implies that
for any $m\le m_0$ we have
$p_m<\eps_d$ (this is because, with probability at least $e^{-\lambda(3m_0)^d}$
 there will be no centers in ${\cal A}(K^m)$, in which case~$K^m$ is not passable). 

Then, if $n>m_0$ and $p_{3\log n}<\eps_d$ we have by~(\ref{pass})
\begin{eqnarray*}
p_n&\le& \Big(\frac{n}{\log n}\Big)^d (11^d p_{3\log n})^{k_0}
  + (3n)^d e^{-\log^d n}\\
&\le& \Big(\frac{n}{\log n}\Big)^d (11^d \eps_d)^{k_0}
  + (3n)^d e^{- \log^d n}\\
&\le& (3n)^d(11^d \eps_d)^{n^{1/2}}+ (3n)^d e^{- \log^d n}\\
&<& \eps_d.
\end{eqnarray*}
By induction, this implies that $p_n<\eps_d$ for all~$n$ (i.e., using the above
calculation, first we obtain that $p_m<\eps_d$ for all $m\le m_0$ implies
that $p_m<\eps_d$ for all $m\leq e^{m_0/3}$, and so on).
Moreover, using~(\ref{pass}) once again, we obtain
\begin{eqnarray}
p_n&\le &(3n)^d(11^d \eps_d)^{n/2}+ (3n)^d e^{- \log^d n}\nonumber\\
&\le &2(3n)^d e^{- \log^d n}, \label{p_n}
\end{eqnarray}
so $p_n\to 0$ as $n\to \infty$.
% So, if $\eps$ is small enough (i.e., $\lambda$ was chosen small), $n$ is large enough,
% and $p_{3\log n}<\eps$, then also
% $p_n<\eps$. 
% 
% For fixed $m_0>1$, for any $\eps>0$ we can choose $\lambda$ small
% enough so that for any $m\le m_0$ we have
% $p_m< \eps$,
% since for~$\lambda$ small with large probability there will
% be no centers in ${\cal A}(K^m_0)$.
Using~(\ref{oc_R1}), one can write (recall that $c(\lambda)>c'(\lambda)\geq 1$)
\begin{eqnarray}
\lefteqn{\Pr[K^n_0 \text{ intersects with a connected component of diameter at 
              least $n/2$}}\nonumber\\
&&~~~~\text{of painted level-$1$ cubes}]\nonumber\\
& \le&
   \Pr[\text{$K^n_0$ is passable}] \nonumber\\
  && {}+\Pr\Big[\text{there exists $i\in{\cal A}(K^n_0)\cap\Z^d$ such that } 
R_i\ge \frac{n}{6(\beta_d+1)}\Big]\nonumber\\
&&  {}+\Pr[\text{there exists $i\in \Z^d\setminus
              {\cal A}(K^n_0)$ such that } R_i\ge \rho(\{i\}, K^n_0)]\nonumber\\
&<& p_n+ (3n)^d e^{-(n/6(\beta_d+1))^d}+c_2 \sum_{\ell=n}^\infty \ell^{d-1}
         \Pr[R_\ell\ge \ell ]\nonumber\\
&\le& p_n+ (3n)^d e^{-(n/6(\beta_d+1))^d}+c_2 \sum_{\ell=n}^\infty
            \ell^{d-1}e^{-\ell^d}\nonumber\\
&\le& 2(3n)^d e^{- \log^d n}+(3n)^d e^{-(n/6(\beta_d+1))^d}+c_2 \sum_{\ell=n}^\infty 
\ell^{d-1}e^{-\ell^d}\label{BC}\\
&\to& 0,\nonumber
\end{eqnarray}
as $n\to \infty$.

We proved that we can choose~$\lambda$ small enough to obtain
\begin{eqnarray*}
&&\Pr[K^n_0 \text{ intersects with a connected component of diameter at 
least $n/2$}\\
&&~~~~\text{of painted level-$1$ cubes}]\to 0,
\end{eqnarray*}
as $n\to \infty$. 
Note that, since ${\cal C} \subset {\mathfrak C}$
\begin{eqnarray*}
&&\Pr[K^n_0 \text{ intersects with a connected component of diameter at 
least $n/2$}\\
&&~~~~\text{of painted level-$1$ cubes}]\\
&\ge& \Pr[0 \text{ belongs to an unbounded connected subset of }{\cal C}]
\end{eqnarray*}
and the latter probability is strictly positive, 
in the case when there is percolation. So, there is no 
percolation for~$\lambda$ small enough and the part~(i) of
Theorem~\ref{t_main} is proved.

\medskip

As for the part~(ii), we proceed as follows. Denote by $H_2\subset \R^d$
the two-dimensional plane:
\[
 H_2 = \{x=(x^{(1)},\ldots,x^{(d)})\in \R^d : x^{(3)}=\ldots =x^{(d)} =0\}.
\]
Now one can write
\begin{eqnarray*}
\lefteqn{\{\text{there is no unbounded connected subset in }\R^d\setminus {\cal C}\}}\\
&&\subset  \{\text{for any bounded~$W\subset H_2$,
there is a contour around $W$ in } {\cal C}\cap H_2\}\\
&&\subset  \{\text{for any bounded~$W\subset H_2$,
there is a contour around $W$ in } {\mathfrak C}\cap H_2\}\\
&&\subset  \{\text{for an infinite number of cubes $K^n_0$, $n=1,2,3,\ldots$,  $K^n_0$ intersects }\\
&&~~~~~ \text{ with a connected component of diameter at least $n/2$}\\
 &&~~~~~ \text{ of painted  level-$1$ cubes}\}.
\end{eqnarray*}
By~(\ref{BC}) and Borel-Cantelli lemma, for small enough~$\lambda$
the probability of the last event is~$0$, and thus part~(ii) of
Theorem~\ref{t_main} is proved.
\qed

\section*{Acknowledgements}
S.P.\ and M.V.\ are grateful to CNPq (302981/02--0 and 306029/03--0)
 for partial support.
M.V.F. acknowledges the support by
CNPq (150989/05--9) and FAPESP (2005/00248--6). This
research was developed using resources of CENAPAD-SP (National Center of High
Performance Computing at S\~{a}o Paulo), project UNICAMP/FINEP-MCT, Brazil.
The authors are grateful to Yuval Peres, who suggested this problem
to them, and to Daniel Andr\'es D\'\i{}az Pach\'on, who has read the first version
of the manuscript very carefully and pointed out several mistakes.
Also, the authors thank the anonymous referee for valuable comments
and suggestions.

\end{document}